%% file: arxiv.tex
\documentclass[letterpaper, 10 pt]{article}

\usepackage{basic_style}
\usepackage{pgfplotstyle}

\usepackage{algorithm,algorithmic}

\newcommand{\tikzsetnextfilename}[1]{%
}
\usepackage{url}

\usepackage{balance}
\usepackage{microtype}

\newcommand{\idxset}{\mathcal{I}}
\newcommand{\Tset}{\mathcal{T}}
\newcommand{\Yset}{\mathcal{Y}}

\newcommand{\onote}[1]{{\color{blue} #1}}

\newcommand{\unote}[1]{}

\title{A Refined Algorithm for Curve Fitting by Segmented Straight Lines
}

\author{Olof Troeng, Mattias Fält%
\thanks{The authors are with the Department of Automatic Control, Lund University, Sweden. E-mail: \texttt{\{oloft,mattiasf\}@control.lth.se}. The authors have received financial support from the Excellence Center at Linköping-Lund in Information Technology (ELLIIT) and the Swedish Research Council through the LCCC Linnaeus Center.}
}

\begin{document}

\maketitle
\thispagestyle{empty}
\pagestyle{empty}

\begin{abstract}

We consider least squares approximation of a function of one variable by a continuous,
piecewise-linear approximand that has a small number of breakpoints.
This problem was notably considered by Bellman who proposed an approximate algorithm based on dynamic programming. Many suboptimal approaches have been suggested, but so far, the only exact methods resort to mixed integer programming with superpolynomial complexity growth.

In this paper, we present an exact and efficient algorithm based on dynamic programming with a hybrid value function. The achieved time-complexity seems to be polynomial.
\end{abstract}

\section{Introduction}
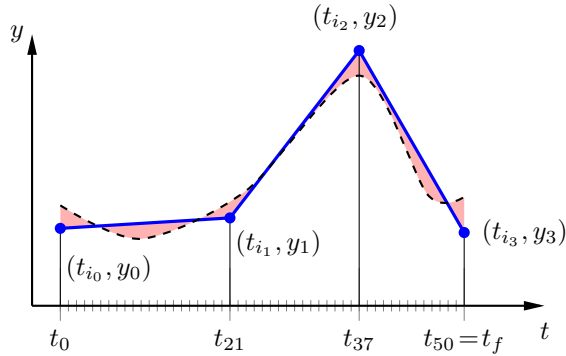
\begin{figure}
	\centering
	\tikzsetnextfilename{img/problem_illustration}%
	\input{img/problem_illustration.tex}%

	\caption{Illustration of the considered problem. For a given function $g$ (dashed line) defined over $[t_0, t_f]$, we want to find a continuous,
	piecewise-linear function $f_{\idxset, \Yset}$ (solid line), with a small number of breakpoints,
	that minimizes the approximation error
	$\norm{ f_{\idxset, \Yset} - g}_2^2 = \int_{t_0}^{t_f} (f_{\idxset, \Yset}(\tau) - g(\tau))^2 d\tau$.
	The breakpoints of $f_{\idxset, \Yset}$ are restricted to belong to some given set $\Tset = \{t_k\}$ (illustrated by the tick marks on the $t$ axis),
	but the $y$ values can be chosen freely. The notation in the figure is further explained in Sec.~\ref{sec:problem-formulation}.}

    \label{fig:approx_err_illustration}
\end{figure}

A classic optimization problem, is that of approximating a function of one variable by a
piecewise-linear approximand with a small number of breakpoints,
that should be selected from a given set (see Fig.~\ref{fig:approx_err_illustration}).
While this problem has been studied extensively,
we are not aware of any previous approach that \emph{efficiently} finds an approximand that is
\emph{optimal} in \emph{2-norm}, \emph{continuous},
and which may take \emph{arbitrary values} at the breakpoints.

The problem was for instance considered  by Bellman \cite{BelRoh69}, who \emph{restricted} the function values of the approximand at the breakpoints to finite sets, enabling a solution by straight-forward dynamic programming.
Similarly, it is easy to solve the problem if the continuity constraint is dropped \cite{Bellman1961}, or if the approximand is restricted to equal the approximated function at the breakpoints \cite{Camponogara2015}.
Another approach is reweighted $\ell_1$ regularization \cite{Candes2008,Kim2009}, which efficiently computes a sparse, but \emph{suboptimal} solution.
Formulating and solving the problem as a mixed-integer quadratic program gives a certifiably optimal solution \cite{Roll2004}, but although the performance of mixed integer solvers have improved significantly over the last decade \cite{Bertsimas2016},
they still have a superpolynomial complexity growth in general.

In this paper, we will present an algorithm for the exact and efficient solution of the considered problem, allowing the approximand to take arbitrary values at the breakpoints.
The algorithm is based on dynamic programming, but more sophisticated than the one in \cite{BelRoh69}---relying on that the hybrid value function can be represented exactly using  piecewise-quadratic functions.
Crucial to the performance, is that the piecewise-quadratic representations can be kept minimal with little computational effort.
The worst-case time complexity seems to have a quartic dependence on the possible number of breakpoints---we have no proof of this,
but extensive numerical tests suggest it.
See \cite{EllZeroTrendFiltering} for an implementation of our algorithm in Julia, that solves medium-size problems
with 1000--10\,000 possible breakpoints in seconds to minutes.

Two notable papers that rely on piecewise-quadratic value functions are:
\cite{Johnson2013}, where they were elegantly used for accelerating $\ell_0$ regularized, piecewise-\emph{constant} regression, and \cite{Lincoln2006}, where \emph{multivariate} piecewise-quadratic function was used for representing the relaxed value function in an optimal control problem.
In the second paper, the S-procedure had to be used to keep the piecewise-quadratic representations small.
In this paper however, where the piecewise quadratics are univariate,
the representations can be kept \emph{minimal} by storing them as linked list.
This is computationally more efficient than relying on the S-procedure.

In Sec.~\ref{sec:approx_algorithm} we present and discuss our algorithm, and in Sec.~\ref{sec:numerical_example} we provide a numerical example.

\emph{Remark:}
The discrete-time (time series) analog of the considered problem is known as $\ell_0$ trend filtering, or segmented-line regression,
and has important applications in a wide range of disciplines---from economy to biology \cite[Sec. 3.1]{Kim2009}.
Our algorithm also handles this discrete problem, but we focus our presentation on the continuous formulation and the related dynamic programming,
since we feel that this is more appealing to a controls audience.
In the seminal paper on $\ell_1$ trend filtering \cite{Kim2009} it was claimed that $\ell_0$ trend filtering was intractable.

\section{Approximation Algorithm}\label{sec:approx_algorithm}
\subsection*{Section Outline}
We begin with a precise formulation of the considered problem  (\mbox{Sec.~\ref{sec:problem-formulation}}). Then, we show that the approximation error between two breakpoints is a quadratic form in the values of the approximand at the breakpoints (Sec.~\ref{sec:linear_fcn_approx_error}). With the help of the derived expression, we can introduce a value function, that once fully computed, gives the solution to the original problem (Sec.~\ref{sec:value_fcn}). Then, we present the following three key points that make the algorithm efficient and exact.
\begin{enumerate}
\item The value function can be represented exactly using a collection of piecewise-quadratic functions (Sec.~\ref{sec:pw_quadratic_representation}).
\item The dynamic programming step of computing one of these piecewise-quadratics from those that have already been computed, is straight forward (Sec.~\ref{sec:dyn_prog_step}).
\item The representations of the piecewise-quadratic functions can be kept minimal with little computational effort (Sec.~\ref{sec:efficient_pwq_representaiton}).
\end{enumerate}
The presentation of the algorithm is concluded by a pseudo-code implementation (Sec.~\ref{sec:algorithm_pseudo_code}).

For the specialized case of (continuous-time) $\ell_0$ trend filtering,
the algorithm can be made significantly more efficient with minor modifications---this we describe in Sec.~\ref{sec:l0_regularized}.
Then, we discuss the question of the computational complexity, and justify that the algorithm is indeed efficient (Sec.~\ref{sec:complexity}).

We conclude the section by briefly introducing the discrete-time analog of the problem (Sec.~\ref{sec:discrete_problem}),
and with the discrete-time formulation we explore the connection to optimal control of hybrid systems (Sec.~\ref{sec:optimal_control_connection}).

\subsection{Problem Formulation}\label{sec:problem-formulation}
Given
\begin{enumerate}
\item a function $g\!: [t_0, t_N] \rightarrow \mathbb{R}$
\item an increasing sequence $\Tset = \{t_0, \ldots, t_N \}$,
\end{enumerate}
we want to find a \emph{short} subsequence $\{t_{i_k}\}_{k=0}^M$ of $\Tset$, indexed by an index set $\idxset := \{i_k\}_0^M \subset \{ 0, \ldots, N\}$,
and a corresponding sequence $\Yset:=\{y_k\}_0^M$ of values,
so that the continuous, piecewise-linear function
\[
f_{\idxset, \Yset}(t) :=
\frac{t_{i_{k+1}}-t}{t_{i_{k+1}} - t_{i_k}} y_k
+
\frac{t-t_{i_k}}{t_{i_{k+1}} - t_{i_k}} y_{k+1},
\quad
t \in [t_{i_k}, t_{i_{k+1}}],
\]
is a good approximation to $g$ in the least squares sense. Note that we implicitly assume that $0$ and $N$ belong to $\idxset$. See Fig.~\ref{fig:approx_err_illustration} for an illustration.

We will consider two versions of this problem, \emph{$\ell_0$ constrained minimization}
\begin{subequations}
    \begin{align}
    \mathop{\text{minimize}}_{\idxset, \Yset} \quad &
    \norm{ f_{\idxset, \Yset} - g}_2^2
    \notag
    \\
    \text{subject to} \quad &
    \text{card}(\idxset) - 1 = M \label{eq:l0-constrained}
    \end{align}
	and \emph{$\ell_0$ regularization}
    \begin{equation}
    \mathop{\text{minimize}}_{\idxset, \Yset} \quad \norm{ f_{\idxset, \Yset} - g}_2^2 + \zeta(\text{card}(\idxset)-1),\label{eq:l0-regularization}
    \end{equation}%
	\label{eq:l0-problem}%
\end{subequations}%
where $\text{card}(S)$ is the number of elements in the set $S$.

We have assumed equality for the cardinality constraint in (\ref{eq:l0-constrained}).
This simplifies the exposition, and can be done without loss of generality since it is clear
that the objective value is non-increasing with $M$.

For the presentation of our algorithm, we will focus on problem (\ref{eq:l0-constrained}),
since this problem is more general, and was considered in \cite{BelRoh69}.

\emph{Remark:} For $\ell_1$-regularization the two formulations
(\ref{eq:l0-constrained}) and (\ref{eq:l0-regularization}) are essentially equivalent:
sweeping $\zeta$ and $M$ generates the same sets of regularizations.
However since the cardinality constraint ($\ell_0$-``norm'') is not convex,
there will typically be solutions to \eqref{eq:l0-constrained} that cannot be reproduced by solving \eqref{eq:l0-regularization}.

\subsection{Approximation Error Between Breakpoints}\label{sec:linear_fcn_approx_error}
The error when the function $g(t)$ is approximated over an interval $[t_i, t_{i'}]$
by an affine function, that takes the values $y$ and $y'$ in the endpoints, is given by
\newcommand{\ci}{\frac{t_{i'}-t}{t_{i'} - t_i}}
\newcommand{\cj}{\frac{t-t_i}{t_{i'} - t_i}}
\newcommand{\integ}[1]{\int_{t_i}^{t_{i'}} #1 \diff t}
\begin{multline}
l_{ii'}(y, y')  :=
\integ{ \left[ g(t) - \left(\ci y + \cj y' \right) \right]^2 }
\\
= \bmat{y & y'}^T P_{ii'} \bmat{y & y'} + q_{ii'}^T \bmat{y & y'} + r_{ii'}.
\label{eq:cont_transition_cost}
\end{multline}
The expressions for $P_{ii'}$, $q_{ii'}$, and $r_{ii'}$ are given in the Appendix. Note that (\ref{eq:cont_transition_cost}) is a quadratic form in $y$ and $y'$.

With this notation we can rewrite (\ref{eq:l0-constrained}) as
\begin{align}
    \mathop{\text{minimize}}_{\idxset, \Yset} \quad &
    \sum_{k=0}^{M-1} l_{i_k i_{k+1}}(y_{k}, y_{k+1})
    \notag
    \\
    \text{subject to} \quad &
    \text{card}(\idxset) -1 = M . \label{eq:l0-constrained_rewrite}
\end{align}

\subsection{Value Function}\label{sec:value_fcn}
We will present an algorithm that relies on dynamic programming to find optimal solutions to \eqref{eq:l0-constrained_rewrite}.

To this end, introduce $V(i,m,y)$ as the minimal cost over the subinterval $[t_i,t_N]$, using $m$ segments,
as a function of the value $y$ at $t_i$, i.e.
\begin{align}
    V(i, m, y) := \quad \mathop{\text{min}}_{\idxset, \Yset} \quad \quad  &
    \sum_{k=0}^{m-1} l_{i_k i_{k+1}}(y_{k}, y_{k+1}) &
    \notag
    \\
     \text{subject to} \quad &
    \text{card}(\idxset) - 1 = m \notag\\
    & \idxset \subset \{i, \ldots, N\} \notag\\
	& y_0 = y, \label{eq:value_fcn_def}
\end{align}
where $\Yset=\{y_k\}_{k=0}^m$, $\idxset =\{i_k\}_{k=0}^m $.

For convenience, we introduce the notation $V_i^m(y) := V(i, m, y)$.
Note that $\min_y V_0^M(y)$ equals the optimal cost of (\ref{eq:l0-constrained_rewrite}).

\subsection{Piecewise-Quadratic Representation of Value Function}\label{sec:pw_quadratic_representation}
\begin{figure}
	\centering
	\tikzsetnextfilename{img/piecewise_quadratic}%
	\input{img/piecewise_quadratic.tex}%

	\caption{The value function $V(m,i,y) = V_i^m(y)$ is piecewise quadratic in the $y$ argument.
	This allows the representation $V^m_i(y) = \min_{p \in \Pi^m_i} p(y)$, where $\Pi^m_i$ is a set of quadratic polynomials.}
	\label{fig:piecewise_quadratic}
\end{figure}
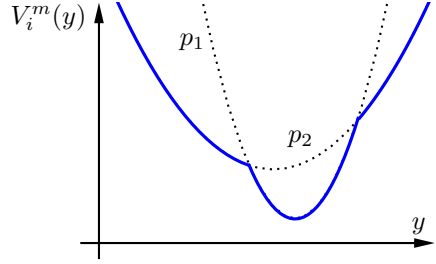

One key reason that our algorithm is able to find the exact solution is that every $V_i^m(y)$ will be piecewise quadratic
\begin{equation}
V^m_i(y) = \min_{p \in \Pi^m_i} p(y),
\label{eq:value_fcn_pwq_representation}
\end{equation}
where the set $\Pi^m_i$ is a set of positive definite, univariate, quadratic functions. See Fig.~\ref{fig:piecewise_quadratic} for an illustration.
That this representation is possible follows from that, for each selection of $\idxset$,
the objective in (\ref{eq:value_fcn_def}) is quadratic in the elements of $\Yset$.
It is clear that minimizing this function with respect to all but the first element in the sequence,
gives a quadratic function.
The function $V^m_i$ is simply the minimum of all such quadratics for different index sets $\idxset$.
This representation will also become clear from the dynamic programming step in the next section.

\subsection{Dynamic Programming Step}\label{sec:dyn_prog_step}
From the definition of the value function (\ref{eq:value_fcn_def}) we get
\begin{subequations}
\begin{align}
V^{m}_i\!(y) & = \hspace{-0.17cm}\min_{i' \in \{i+1,..,N-m\}} \, \min_{y'} \, l_{ii'}\left( y , \, y' \right)
+ V^{m-1}_{i'}(y')  \label{eq:dynprog_step_ver1}\\
\intertext{for $m\geq2$, with the initial case for $m=1$ segment:}
V^1_i(y) & = \min_{y'} \, \left[ l_{iN}(y, y') + V^0_N(y')\right], \quad V_N^0(y)=0 \label{eq:dynprog_base}
\end{align}
\end{subequations}
Although the terminal cost $V_N^0$ is identically zero,
we have included it for consistency with the dynamic programming literature,
and to allow for a natural transition to a discrete formulation of the optimization problem \eqref{eq:l0-problem}.

\begin{figure}
	\centering
	\vspace{0.8em}
	\resizebox{\columnwidth}{!}{
	\tikzsetnextfilename{img/multistage_graph}%
	\input{img/multistage_graph.tex}%

	}
	\caption{Illustration of the dynamic programming problem as a multistage graph when the set $\Tset=\{t_k\}_{k=0}^5$
	contains 6 possible breakpoints,
	of which 2 are to be selected ($t_0$ and $t_N$ are always included).
	Note that the cost $V^m_i$ at each node is given by a piecewise-quadratic function and that the transition costs
	$l_{ii'}$ are quadratic forms.
	Following the tradition of optimal control, the search is made backwards in time, starting from $V_5^0$.}
	\label{fig:multistage_graph}
\end{figure}
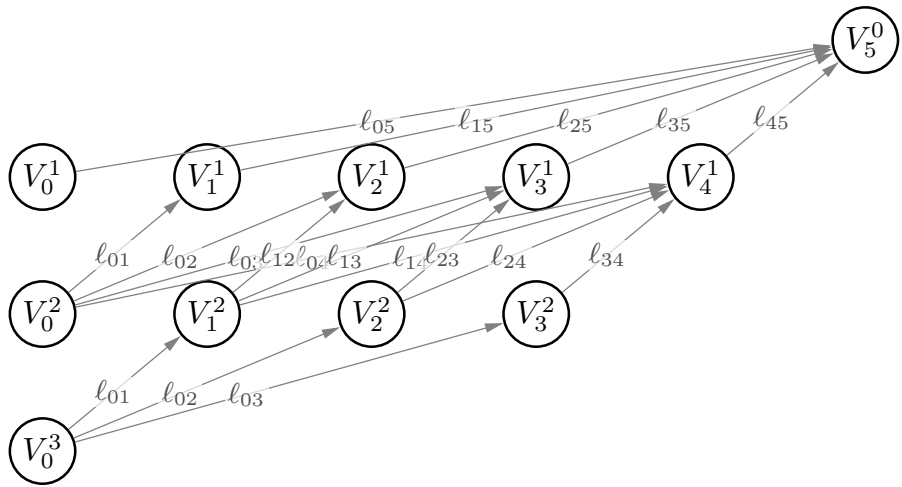

It is now clear that we, in principle, could compute $V^{m+1}_i$
from knowledge of $V^{m}_{i'}$, $i' > i$ in a dynamic programming fashion,
according to Fig.~\ref{fig:multistage_graph}.
Next, we show that the minimization in (\ref{eq:dynprog_step_ver1}) actually is easy to perform.

Using the representation (\ref{eq:value_fcn_pwq_representation}) we can write (\ref{eq:dynprog_step_ver1}) as
\begin{align}
	V^{m+1}_i(y) & = \hspace{-0.17cm}\min_{i' \in \{i+1,..,N\!-m\}} \,\, \min_{y'} \, \left[ l_{ii'} ( y , \, y' ) + V^{m}_{i'}(y') \right]
\notag \\
& = \hspace{-0.17cm}\min_{i' \in \{i+1,..,N\!-m\}} \,\, \min_{y'} \, \min_{p \in \Pi^{m}_{i'}} \, \left[ l_{ii'}( y , \, y' ) +  p(y') \right]
\notag  \\
& = \hspace{-0.17cm}\min_{i' \in \{i+1,..,N\!-m\}} \, \min_{p \in \Pi^{m}_{i'}} \, \left[ \min_{y'}  \, \left[ l_{ii'}( y , \, y') +  p(y') \right] \right]
\label{eq:dynprog_step_ver2}
\end{align}
The expression in the innermost parentheses is a positive definite quadratic form in $y$ and $y'$,
so the minimum with respect to $y'$ is a quadratic function that is easy to compute exactly.
Relation \eqref{eq:dynprog_step_ver2} shows that the representation $\Pi^{m+1}_i$ of $V^{m+1}_i$
in \eqref{eq:value_fcn_pwq_representation}, is straight forward to compute from previous $\Pi^m_{i'}$.
The possibility of swapping the min operations in the last equality above,
enables us to work with $\Pi^{m+1}_{i}$ as the representation of $V^{m+1}_i$---this is a key step to make the algorithm efficient.

\subsection{Minimal Representation $\Pi_i^m$ of $V_i^m$}\label{sec:efficient_pwq_representaiton}
\emph{In principle} $\Pi^{m+1}_i$ could be generated directly from $\left\{ \Pi^m_{i'} \right\}_{i' \in \{i+1,..,N\!-m\}}$ as
\begin{equation} \Pi^{m+1}_i =\hspace{-0.17cm}\bigcup_{i' \in \{i+1,..,N\!-m\}} \, \bigcup_{p \in \Pi^{m}_{i'}}
\left\{   \min_{y'}  \, \left[ l_{ii'}(y, \, y') + p(y') \right] \right\}.
\label{eq:union_of_polynomials}
\end{equation}
If (\ref{eq:union_of_polynomials}) would simply be iterated, then $\Pi^m_i$ would grow rapidly in size.
So the third, crucial step for an efficient implementation is to keep the size of $\Pi^m_i$ as small as possible.
Since $V_i^m(y)=\min_{p\in\Pi_i^m}\,p(y)$, it is sufficient to keep only those polynomials
that are smaller than all other polynomials in $\Pi^m_i$ over some interval.

\begin{figure}
    \centering
    \tikzsetnextfilename{img/pwp_illustration1}
    \input{img/pwp_illustration}
	\caption{
	Representing the piecewise-quadratic function $V_i^m(y)$ (blue) by a linked list allows
	for efficient insertion of new polynomials.
	Each list node contains an interval
	and the quadratic polynomial that attains the minimum over this interval.
	In this example, the list would contain the nodes
	$([-\infty,\alpha_1],p_2), \,([\alpha_1,\alpha_2],p_1),$ and $([\alpha_2,\infty],p_2)$.
	To add a new quadratic $\mu$ to the representation,
	it is sufficient to look for intersections with each of the polynomials on their corresponding intervals.
	In this example, $p_1$ would simply be replaced with $\mu$ and the intervals updated.
	}
    \label{fig:pwp_illustration}
\end{figure}
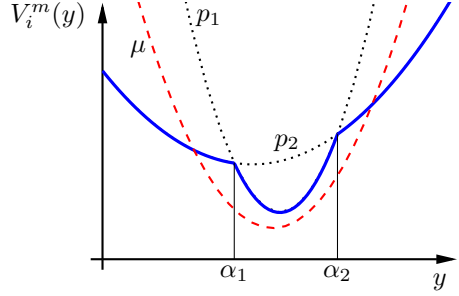

This can be implemented efficiently using linked lists,
where each list node keeps track of an interval and a quadratic polynomial,
see Fig.~\ref{fig:pwp_illustration} for an illustration.
With slight abuse of notation, we will refer to both the sets and linked lists using the notation $\Pi_i^m$.
Starting with an empty linked list $\Pi^{m+1}_i$, the elements in \eqref{eq:union_of_polynomials} are inserted one by one,
while redundant list nodes are discarded.
Compare this to \cite{Lincoln2006}, where the S-procedure was used to obtain parsimonious, but typically non-minimal, representations.

\subsection{Algorithm}\label{sec:algorithm_pseudo_code}
The algorithm resulting from the discussion above is described in pseudo-code in Algorithm~\ref{alg:optimal_fit}.

Recovering the solution (line \ref{alg:constrained-recover})
is a standard step in dynamic programming, and is trivial to do if each
polynomial $p$ keeps track of the index $i'$ when it was generated.

Note that the algorithm computes $V_0^m$ for all $m\leq M$,
i.e., we get all solutions to problem~\eqref{eq:l0-constrained}
with less than $M$ segments for free.

\renewcommand{\algorithmicrequire}{\textbf{Input:}}
\renewcommand{\algorithmicensure}{\textbf{Output:}}
\begin{algorithm}[H]
	\caption{Continuous Piecewise-Linear Approximation}
	\label{alg:constrained}
	\begin{algorithmic}[1]
		\REQUIRE \parbox[t]{6cm}{
			Function $g(\cdot)$ to approximate\\
			Time points $\Tset = \{t_0, t_1, \ldots, t_N\}$
		}
		\vspace{0.2em}
		\ENSURE Solution $(\idxset$, $\Yset)$ of (\ref{eq:l0-constrained})
		\vspace{0.2em}

		\STATE Compute transition costs $l_{ii'}( \cdot , \, \cdot )$  from $\Tset$ and $g(\cdot)$\\
		\quad \emph{($l_{ii'}(\cdot,\cdot)$ are quadratic forms)}

		\vspace{0.2em}
		\FOR{$i = 0 \textbf{ to } N-1$}
		\STATE Let  $\Pi^{i}_{N} = \left\{ \min_{y'} l_{iN}( \cdot , \, y' ) \right\}$
			\emph{(representation of $V^{i}_N$)}
		\ENDFOR

		\vspace{0.5em}
		\FOR{$m = 2 \textbf{ to } M$}
		\FOR{$i = 0 \textbf{ to } N-m$}
		\STATE Initialize $\Pi^{m}_{i} = \{\}$
		\FOR{$i' = i+1 \textbf{ to } N-m+1$}

			\FOR{$p \textbf{ in } \Pi^{m-1}_{i'}$}
			\STATE Compute $\mu(\cdot) = \min_{y'}  \, l_{ii'}( \cdot , \, y'  ) + p(y')$\\
			\quad \emph{(It is clear that $\mu$ is a quadratic)}
			\vspace{0.2em}
			\STATE Add $\mu$ to  $\Pi_i^{m}$ \label{alg:constrained-addquadratic}\\
			\quad \emph{(Keep only if $\mu$ is smallest at some interval)}
			\ENDFOR
		\ENDFOR
		\\\emph{(We now have a minimal representation $\Pi_i^{m}$ of $V_i^{m}$)}
		\ENDFOR
		\ENDFOR
	\STATE Find $\text{argmin}_{p\in \Pi^M_0} \min_{y'}  \, p( y' )$
	\\\quad\emph{(i.e $p$ corresponding to minimum of $ V^{M}_{0}$)}
	\STATE Recover the corresponding solution $\idxset$, $\Yset$ \label{alg:constrained-recover}
	\end{algorithmic}

    \label{alg:optimal_fit}
\end{algorithm}

\subsection{$\ell_0$ Regularization}\label{sec:l0_regularized}

Although Algorithm~\ref{alg:optimal_fit} can be used for also solving \eqref{eq:l0-regularization}, the efficiency can be improved by specializing the algorithm.

This is done by instead considering the following value function without an explicit dependence on $m$
\begin{align*}
V(i, y) := \quad \mathop{\text{min}}_{\idxset, \Yset} \quad  &
\hspace{-0.3cm}\sum_{k=0}^{\text{card}(\idxset)-2}\!\!\! l_{i_k i_{k+1}}(y_{k}, y_{k+1}) + \zeta (\text{card}(\idxset)-1)&
\notag
\\
\text{subject to} \quad & \idxset \subset \{i, \ldots, N\} \notag\\
& y_0 = y. \label{eq:value_fcn_def}
\end{align*}
Introduce the notation $V_i(k) := V(i, y)$, the dynamic programming step takes the form
\begin{align*}
V_i(y) & = \hspace{-0.17cm}\min_{i' \in \{i+1,..,N\}} \,\, \min_{y'} \, \left[ l_{ii'} ( y , \, y' ) + V_{i'}(y') + \zeta \right]
\notag \\
& = \hspace{-0.17cm}\min_{i' \in \{i+1,..,N\}} \, \min_{p \in \Pi_{i'}} \, \left[ \min_{y'}  \, \left[ l_{ii'}( y , \, y') +  p(y') + \zeta \right] \right],
\end{align*}
where $V_N(y)=0$.

The main difference to the constrained problem,
is that there is no longer an explicit dependence on the number of segments $m$,
this is instead captured by the cost $\zeta$.
See the implementation \cite{EllZeroTrendFiltering} for further details.
\unote{Remove explicit defenition}

\subsection{Complexity}\label{sec:complexity}

The representation of the piecewise-quadratics are linked lists where each element
is a quadratic function, and the interval over which this function defines the piecewise-quadratic.
Because of this representation, when adding a new quadratic function,
it is easy to go through the list, find potential intersections on the respective intervals,
and add the new quadratic function if needed.
The complexity for insertion into the linked list $\Pi_i^m$
(Algorithm \ref{alg:constrained}, line \ref{alg:constrained-addquadratic})
is thus linear in the list length.

\begin{figure}
    \centering
    \newcommand{\dpath}{img}
	\tikzsetnextfilename{img/complexity}%
	\input{img/complexity.tex}%

    \caption{Illustration of the number of segments in the piecewise-quadratics
		for different problems. Each gray line shows $\max_m\text{length}(\Pi^m_i)$,
		as a function of $i$, for a specific problem.
		The problem data were generated from polynomials, exponentials,
		random sequences, as well as synthetically generated with circular arcs and sharp kinks.
		The red (dashed) line represents the apparent bound $\text{length}(\Pi^m_i) \leq N-i$.
		}
    \label{fig:nbr_segments}
\end{figure}
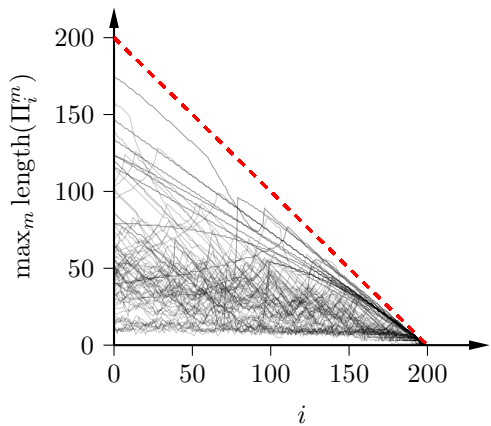

If we let $R$ be the maximum length over all $\Pi_i^{m}$, at \emph{any} time in the algorithm,
we see from Algorithm~\ref{alg:constrained} that we get the following bound on the complexity
\newcommand{\Oh}{\mathcal{O}}
\[
\mathcal{O} (MN^2R^2).
\]
The questions is if it is possible to bound $R$? Although we have no proof, extensive testing suggests that $R \leq N$.
A plot of the list lengths, for a wide class of problems, is shown in Fig.~\ref{fig:nbr_segments}.
Assuming that the conjectured inequality holds, the algorithm for $\ell_0$ constrained approximation has a worst-case time complexity of $\mathcal{O} (MN^4)$.

In the algorithm for $\ell_0$ regularization problem, there is no outer loop over $m$.
This gives the complexity $\mathcal{O}(N^2R^2)$.

\subsection{The Discrete Problem}\label{sec:discrete_problem}
The algorithm can equally well find the optimal approximation to a \emph{time series} $g \in \mathbb{R}^{N+1}$, by an approximand
\[
f_{\idxset, \Yset}[t] :=
\frac{i_{k+1} - t}{i_{k+1} - i_k} y_k
+
\frac{t-i_k}{i_{k+1} - i_k} y_{k+1},
\quad
i_k \leq t \leq i_{k+1}.
\]
The optimization problem still takes form \eqref{eq:l0-constrained} or \eqref{eq:l0-regularization}.
Since the possible breakpoints are simply given by the set $\Tset = \{0, \ldots, N\}$,
we are able to write $i_k$ instead of $t_{i_k}$ in the expression above.

Just as before in \eqref{eq:cont_transition_cost},
the approximation error between two breakpoints $i$ and $i'$ of $f_{\idxset, \Yset}$,
is a quadratic form
\newcommand{\sumexpr}{\sum_{k=i}^{i'-1}}
\newcommand{\Ki}{\frac{i'-k}{i'-i}}
\newcommand{\Kip}{\frac{k-i}{i'-i}}
\begin{align*}
l_{ii'}(y, y') & =
\sumexpr \left[ g[t] - \left(\Ki y + \Kip y' \right) \right]^2
\\
& = \bmat{y & y'}^T P_{ii'} \bmat{y & y'} + q_{ii'}^T \bmat{y & y'} + r_{ii'},
\end{align*}
where
$P_{ii'}$, $q_{ii'}$, and $r_{ii'}$ are given in the Appendix.
The sum above does not include $i'$, since the cost for that index is handled by the succeeding segment.

The only other difference is that we in \eqref{eq:dynprog_base} have
\[V_N^0(y) = (g[N] - y)^2.\]

\subsection{Connection to Optimal Control}\label{sec:optimal_control_connection}
The formulation in the previous section corresponds to optimal control of a discrete double integrator,
with the control signal constrained to $M-1$ impulses,
and the objective of tracking a reference $g[k]$.
This is captured by the cost functional
\begin{equation*}
J(u) = \sum_{i=0}^N (g[k] - x_1[k])^2,
\end{equation*}
the dynamics
\begin{align*}
x_1[k+1] &= x_1[k] + x_2[k] \\
x_2[k+1] &= \qquad \quad \,\, x_2[k] + u[k] \\[0.5em]
m[k+1] &= m[k] + \begin{cases}1 & \text{if } u[k] \neq 0 \\ 0  & \text{if } u[k] = 0 \end{cases},
\end{align*}
the initial conditions
\[
m[0] = 0, \qquad x_1[0], x_2[0] \text{ free},
\]
the final time $N$, and the final set
\[
\Gamma = \mathbb{R} \times \mathbb{R} \times \{M-1\}.
\]

From this we see that the considered problem is related to the rich literature on optimal control of hybrid systems \cite{Hedlund1999,Borrelli2003},
and poses the question if similar problems can be solved efficiently.

\section{Numerical Example} \label{sec:numerical_example}
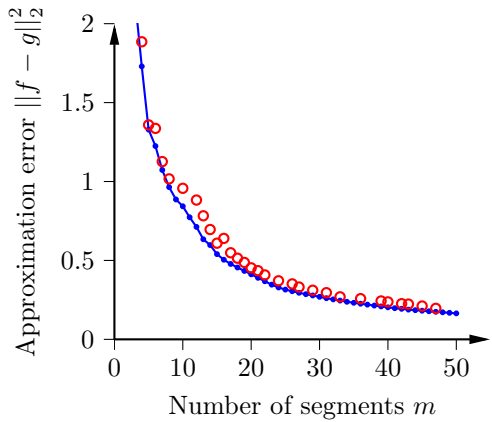
\begin{figure}
	\vspace{0.25em}
	\centering
	\newcommand{\dpath}{img/snp_problem}%
	\tikzsetnextfilename{\dpath/l0_vs_l1}%
	\input{\dpath/l0_vs_l1.tex}%

	\caption{Approximation error vs. the number of segments used, for piecewise-linear approximation
 		of the time series in Sec.~\ref{sec:numerical_example}.
		The blue line correspond to optimal approximands obtained by the discussed algorithm, and the red circles correspond to approximands that were obtained through reweighted $\ell_1$ trend filtering. For a given value of $m$, the optimal approximand has about \SI{15}{\percent} lower error.}
	\label{fig:snp-problem_cost_comparison}
\end{figure}

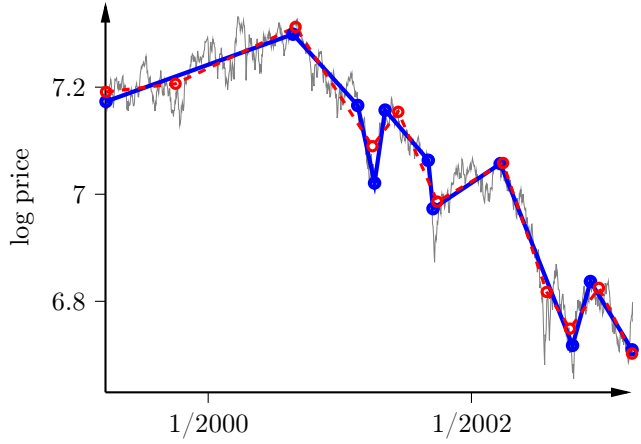
\begin{figure}
	\vspace{0.25em}
	\centering
	\newcommand{\dpath}{img/snp_problem}%
	\tikzsetnextfilename{\dpath/data_fit}%
	\input{\dpath/data_fit.tex}%

	\caption{A time series of 1000 closing prices for the stock index S\&P\,500, together with two piecewise-linear approximands with 10 segments.
		The optimal solution (blue) achieves an error of 0.84, and the approximand obtained through reweighted $\ell_1$ trend filtering (dashed red) achieves an error of 0.96.}
	\label{fig:snp-problem_time_series}
\end{figure}

We have compared the optimal, sparse approximands from our
algorithm to those generated by the popular method of reweighed $\ell_1$ trend filtering
\cite{Candes2008} \cite[Sec. 7.2]{Kim2009}.\footnote{The results in Fig. \ref{fig:snp-problem_time_series} and \ref{fig:snp-problem_cost_comparison}
were generated by sweeping the parameter $s$ in \cite[Sec. 7.2]{Kim2009},
and for each value of $s$ perform 12 iterations for each $\epsilon \in \{10^{-5},3\cdot 10^{-5}\}$, and then selecting the best result.}

We considered the same data as in \cite[Sec. 4]{Kim2009},
which is a time-series with 2000 consecutive closing values for the American stock index S\&P\,500.
The objective was to minimize the approximation error, as a function of the number of segments of the approximand.
\unote{något om log av points, något om fel datum?}

We were not able to make the reweighted $\ell_1$ trend filtering robust enough to handle the full data set.
Therefore we restricted the problem to the first 1000 data points.

The cost vs. the number of segments in the piecewise-linear approximation is shown in Fig.~\ref{fig:snp-problem_cost_comparison}.
The data, together with the optimal and $\ell_1$ approximands are shown in Fig.~\ref{fig:snp-problem_time_series}
for the solutions with $10$ segments.

For the \emph{full} dataset, it took about 5 minutes for a midrange laptop
to compute the optimal solutions
for \emph{all} values $2 \leq M \leq 50$, using our algorithm for $\ell_0$ constrained minimization.
For the case of $\ell_0$ regularization,
it took between \SI{8}{s} for $\zeta = 0.2 \leftrightarrow \text{card}\,\idxset = 8$,
to \SI{1.0}{s} for $\zeta = 0.01 \leftrightarrow \text{card}\,\idxset = 39$.

\section{Discussion and Outlook}
We have presented an exact and efficient algorithm for continuous, piecewise-linear approximation, with a minimal number of breakpoints. It is interesting to note that  $\ell_0$ regularization was claimed to be intractable in \cite{Kim2009}, but that our algorithm solves medium-sized problems of this type in seconds. Despite this, it is an open questions whether the algorithm has a guaranteed polynomial time complexity---our investigations this far, indicate that it does.

Our algorithm enables comparisons between $\ell_0$ and $\ell_1$ regularization for larger problems than those considered in \cite{Hastie2017},
(in the specialized case of trend filtering).

Two possible extensions of the algorithm are: (1) to optimal, sparse, spline approximation by implementing the S-procedure as in \cite{Lincoln2006}---although the computations would be relatively demanding; (2) to handle more general linear systems in the optimal-control setting of Sec.~\ref{sec:optimal_control_connection}.

\section*{Appendix}

\subsection{Expressions for Continuous-Time Transition Costs}
\label{sec:cont_transition_costs_coeffs}

The coefficients in the transitions costs in (\ref{eq:cont_transition_cost}) are given by
\begin{subequations}
	\begin{align*}
	P_{ii'} & =
	\bmat{
		\integ{  \left(\ci\right)^2  } & \integ{  \ci \cj } \\
		\integ{  \ci \cj  } & \integ{\left(\cj\right)^2  } }
	\notag
	\\[0.1em]
	& \qquad\qquad\qquad\qquad\qquad\qquad=
	\frac{(t_{i'} - t_i)}{6}
	\bmat{2 & 1 \\
		1 & 2 },
	\\[0.3em]
	q_{ii'} & =
	-2\bmat{\integ{  \ci g(t) } \\  \integ{  \cj g(t)  } }
	, \quad
	r_{ii'} =
	\integ{  g(t)^2 } .
	\end{align*}
	\label{eq:continuous_transition_cost_qf}
\end{subequations}

By storing the $3N$ values $H_j[i] = \sum_{k=0}^i h_j[k]$,
where $h_1[i]=\int^{t_{i+1}}_{t_i} g(t) dt$, $h_2[i]=\int^{t_{i+1}}_{t_i} tg(t) dt$,
and $h_3[i]=\int^{t_{i+1}}_{t_i} g(t)^2 dt$ for $0 \leq i < N$,
all $q_{ii'}$ and $r_{ii'}$ can be computed from from $H_j[i']-H_j[i]$ in $\mathcal{O}(1)$.

\subsection{Expressions for Discrete-Time Transition Costs}\label{sec:discrete_transition_costs_coeffs}
The coefficients in the transitions costs in Sec.~\ref{sec:discrete_problem} are given by
\begin{subequations}
	\begin{align*}
	P_{ii'} & =
	\bmat{\sumexpr  \left(\Ki\right)^2   & \sumexpr  \Ki \Kip  \\
		\sumexpr  \Ki \Kip  & \sumexpr  \left(\Kip\right)^2   },
	\\[0.3em]
	q_{ii'} & =
	-2\bmat{\sumexpr \Ki g[t] \\  \sumexpr \Kip g[t] }
	, \quad
	r_{ii'} =
	\sumexpr g[t]^2 .
	\end{align*}
	\label{eq:discrete_transition_cost_qf}
\end{subequations}

These can be computed efficiently, as in the continuous case,
by replacing the integrals by appropriate sums.

\section*{Acknowledgment}
The authors thank their colleagues Bo Bernhardsson, Pontus Giselsson, Anders Rantzer, and Fredrik Bagge Carlson, for helpful comments and suggestions.

\input{arxiv.bbl}
\end{document}

%% file: img/problem_illustration.tex
	\pgfplotsset{compat=1.13}
	\begin{tikzpicture}
	\begin{axis}[
	axis x line=center,
	axis y line=left,
	width=8.4cm,
	height=5.2cm,
	ymin=0,
	ymax=1.3,
	xmin=-0.07,
	xmax=1.2,
	every axis y label/.style={at={(current axis.north west)},xshift=-2mm},
	xlabel style={below,yshift=-1mm},
	xlabel=$t$,
	ylabel={$y$},
	axis line style={thick,->},
    xtick={0,0.02,...,1},
	extra x ticks={0,0.42,0.74,1.0},
    extra x tick labels={$t_0$, $t_{21}$, $t_{37}$, $t_{50}\!=\!t_f$},
    extra tick style={
       	major tick length=7mm,
       	yshift=1.5mm,
       	very thick,
       	color=black
    },
	xticklabels={},
	ytick=\empty,
	clip=false	]

	\newcommand{\xcoords}{(0, 0.37) (0.42, 0.42) (0.74, 1.22) (1, 0.35)}
	\addplot[name path=A,black,very thick,mark=*, color=blue,mark size=1.5pt, mark options={solid}] coordinates
	{\xcoords};

	\addplot[ycomb] coordinates {\xcoords};

	\addplot[smooth, name path=B,black,thick, dashed] coordinates
	{(0,0.48) (0.20,0.32) (0.46, 0.55) (0.74, 1.10) (0.91, 0.53) (1.0, 0.52)};	

	\addplot[pattern=vertical lines, color=red!30] fill between[of=A and B];	

	\node at (0.12, 0.18) {$(t_{i_0}, y_0)$};
	\node at (0.54, 0.3) {$(t_{i_1}, y_1)$};
	\node at (0.73, 1.38) {$(t_{i_2}, y_2)$};
	\node at (1.15, 0.35) {$(t_{i_3}, y_3)$};

	\end{axis}
\end{tikzpicture}

%% file: img/piecewise_quadratic.tex
\begin{tikzpicture}
	\pgfplotsset{compat=1.13}
	\newcommand{\shift}{11} 
	\newcommand{\xmin}{-12} 
	\newcommand{\xmax}{7} 

	\begin{axis}[
	axis x line=center, 
	axis y line=center,
	axis line style={thick,->},	
	height=5cm,
	width=6.3cm,
	xmin=\xmin+\shift,
	xmax=\xmax+\shift,
	ymax=30,
	ymin=-2,
	xlabel=$y$, 
	ylabel={$V^m_i(y)$},
	ylabel style={yshift=1mm, xshift=-13mm},
	xtick=\empty, 
	ytick=\empty,
	clip mode=individual	
	]

	\newcommand{\polya}{0.3*(x-\shift)^2+1*(x-\shift)+10}
	\newcommand{\polyb}{1.1*(x+0.5-\shift)^2 + 3}	
	
	\newcommand{\xa}{\xmin+\shift}		
	\newcommand{\xb}{-3.0+\shift}			
	\newcommand{\xc}{2.9+\shift}
	\newcommand{\xd}{\xmax++\shift}

	\foreach \poly in 
	{\polya, \polyb }{
		\addplot[domain=\xmin+\shift:\xmax+\shift, dotted, thick] {\poly};
	}		
		
	\foreach \poly / \le / \re in 
	{\polya /				\xa / \xb,
	 \polyb	/				\xb / \xc,
	 \polya	/			 	\xc / \xd
	 }  
	{

	\addplot[domain=\le:\re, blue, very thick] {\poly};		
	}
	
	\node at (-6 + \shift, 25) {$p_1$};
	\node at (-0.2 + \shift, 13) {$p_2$};

\end{axis}
\end{tikzpicture}

%% file: img/multistage_graph.tex
\begin{tikzpicture}
\tikzset{
	<-/.style={-{Triangle[length=2mm,width=1mm]}}
}

\pgfmathtruncatemacro{\M}{3}
\pgfmathtruncatemacro{\N}{5}

\providecommand{\hlone}{3}
\providecommand{\hltwo}{1}

\newcommand{\col}{black}
\newcommand{\linestyle}{solid}
\newcommand{\textcol}{black}

\foreach \m in {0,...,\M}
{

\pgfmathtruncatemacro{\iMin}{0}
\pgfmathtruncatemacro{\iMax}{\N-\m}

\ifthenelse{\m=0}{\pgfmathtruncatemacro{\iMin}{\N}}{}
\ifthenelse{\m=\M}{\pgfmathtruncatemacro{\iMax}{0}}{}

\foreach \i in {\iMin,...,\iMax}
{
\node[draw,\linestyle] (\m-\i) [circle,inner
sep=1.4pt,thick,\col] at (\i*1.8, -\m*1.5) {$V^{\m}_\i$};
}

}

\foreach \m in {0,1,2}
{

	\pgfmathtruncatemacro{\mnext}{\m + 1}
	\pgfmathtruncatemacro{\iMin}{0}
	\pgfmathtruncatemacro{\iMax}{\N-\mnext}

	\ifthenelse{\mnext=\M}{\pgfmathtruncatemacro{\iMax}{0}}{}

    \foreach \i in {\iMin,...,\iMax}
    {

        \pgfmathtruncatemacro{\ipMin}{\i+1}
        \pgfmathtruncatemacro{\ipMax}{\N-\m}
        \ifthenelse{\m=0}{\pgfmathtruncatemacro{\ipMin}{\N}}{}

        \foreach \ip in {\ipMin,...,\ipMax}
        {

        \draw[<-,solid, gray] (\mnext-\i) --
        node[pos=0.4,color=black,fill=white,fill opacity=0.7,text
        opacity=1,font=\small,inner sep=1pt]{$\ell_{\i\ip}$}
        (\m-\ip);

        }
    }

}

\ifdefined\usehighlight

\foreach \m/\i in {0/\N, 1/\hlone, 2/\hltwo, 3/0}
{
		\pgfmathtruncatemacro{\mlabel}{\m}
		\node[draw,\linestyle] (\m-\i) [circle,inner
		sep=1.4pt,ultra thick,blue] at (\i*1.8, -\m*1.5) {$V^{\m}_\i$};
}

\foreach \m/\i/\ip in {0/\hlone/\N, 1/\hltwo/\hlone, 2/0/\hltwo}
{
	\pgfmathtruncatemacro{\mnext}{\m + 1}
	\draw[<-,ultra thick, blue] (\mnext-\i) --
	node[pos=0.4,color=blue,fill=white,fill opacity=0.7,text
	opacity=1,font=\small,inner sep=1.5pt] {$\ell_{\i\ip}$}
	 (\m-\ip);
}

\fi

\end{tikzpicture}

%% file: img/pwp_illustration.tex
\begin{tikzpicture}
	\pgfplotsset{compat=1.13}
	\newcommand{\xmin}{-8} 
	\newcommand{\xmax}{7} 
	\newcommand{\ymin}{-3} 
	\newcommand{\ymax}{30} 

	\begin{axis}[
	axis x line=left,
	axis y line=left,
	axis line style={thick,->},
	height=5cm,
	width=6.3cm,
	xmin=\xmin,
	xmax=\xmax,
	ymax=\ymax,
	ymin=\ymin,
	xlabel=$y$,
	ylabel={$V^m_i(y)$},
	ylabel style={at={(-0.02,0.95)}, rotate=-90},
	xlabel style={at={(0.95,-0.02)}},
	xtick=\empty, 
	ytick=\empty,
	clip mode=individual
	]

		\newcommand{\polya}{0.3*(x)^2+1*(x)+10}
		\newcommand{\polyc}{1.7*(x+0.5)^2 + 3}

		\newcommand{\polyd}{0.9*(x+0.8)^2 + 1}

		\newcommand{\xa}{\xmin}
		\newcommand{\xd}{\xmax}

		\newcommand{\xb}{-2.4315}
		\newcommand{\xc}{1.9315}

		\foreach \poly in
		{\polya, \polyc }{
			\addplot[domain=\xmin:\xmax, dotted, thick] {\poly};
		}

		\foreach \poly / \le / \re in
		{\polya /				\xa / \xb,
		 \polyc	/				\xb / \xc,
		 \polya	/			 	\xc / \xd
		 }
		{

		\addplot[domain=\le:\re, blue, very thick] {\poly};
		}

		\draw[] (axis cs:\xb,\ymin) -- (axis cs:\xb,9.342157675);
		\draw[] (axis cs:\xc,\ymin) -- (axis cs:\xc,13.050707675);

		\coordinate (label1) at (\xb,\ymin - 1.8);
		\coordinate (label2) at (\xc,\ymin - 1.8);

		\node[thick] at (label1) {$\alpha_1$};
		\node[thick] at (label2) {$\alpha_2$};

		\node at (-3.5, 28) {$p_1$};
		\node at (-0.2, 12) {$p_2$};

		\node at (-6.5, 24) {$\mu$};

		\addplot[domain=-6.5:5.0, red, dashed, thick] {\polyd};

		\draw[thick] (\xmin-0.5,\ymin) -- (\xmin,\ymin);
		\draw[thick] (\xmin,\ymin-1.35) -- (\xmin,\ymin);
	\end{axis}
\end{tikzpicture}

%% file: img/complexity.tex
\begin{tikzpicture}
\begin{axis}[%
arrowplotLongLabel,  
width=5cm,
height=4.5cm,
scale only axis,
xmin=0,
xmax=240,
ymin=0,
ymax=220,
xlabel={$i$},
xlabel style={at={(0.5,-0.15)}, anchor=north},
ylabel={$\max_m \text{length}(\Pi_i^m$)} ,
]

\foreach \i in {0,...,91}
{
\addplot [black, thin, opacity=0.2]
table[x expr=\coordindex, y expr=\thisrowno{\i}, col sep=comma]{\dpath/complexity.csv};

\addplot [dashed, red, domain=0:200] {200-x};
}

\end{axis}
\end{tikzpicture}

%% file: img/snp_problem/l0_vs_l1.tex
\begin{tikzpicture}
\begin{axis}[%
arrowplotLongLabel,  
width=5cm,
height=4.2cm,
scale only axis,
xmin=0,
xmax=55,
ymin=0,
ymax=2,
xlabel={Number of segments $m$},
xlabel style={at={(0.5,-0.15)}, anchor=north},
ylabel={Approximation error $\norm{f - g}^2_2$},
]

\addplot [thick, blue, mark=*, mark size=0.7]
table[x=m, y=cost, col sep=semicolon] {\dpath/cost_vs_m_l0.csv};

\addplot [thick, red, mark=o, only marks, mark options={solid}]
table[x=m, y=cost, col sep=semicolon, mark=*] {\dpath/cost_vs_m_l1.csv};

\end{axis}
\end{tikzpicture}

%% file: img/snp_problem/data_fit.tex
\begin{tikzpicture}
\begin{axis}[%
arrowplotLongLabel,  
width=7cm,
height=5.2cm,
scale only axis,
xmin=0,
xmax=1000,
ymin=6.63,
ymax=7.36,
xtick={195,695,1199,1703}, 
xticklabels={1/2000, 1/2002, 1/2004, 1/2006},
ylabel={log price},
]
\addplot [black, thin, opacity=0.5]
table[x expr=\coordindex, y index=0] {\dpath/snp500.txt};

\addplot [blue, ultra thick, mark=o]
table[x=I, y=Y, col sep=semicolon] {\dpath/sol_10_l0.csv};

\addplot [red, very thick, dashed, mark=o, mark options={solid}]
table[x=I, y=Y, col sep=semicolon] {\dpath/sol_10_l1.csv};

\end{axis}
\end{tikzpicture}